\documentclass[article, 12pt]{imsart}
\usepackage{graphicx}
\graphicspath{ {images/} }
\RequirePackage[OT1]{fontenc}
\RequirePackage{amsthm,amsmath, amssymb, enumerate}
\RequirePackage[numbers]{natbib}
\RequirePackage[colorlinks,citecolor=blue,urlcolor=blue]{hyperref}
\usepackage[english]{babel}

\startlocaldefs
\numberwithin{equation}{section}
\theoremstyle{plain}

\endlocaldefs

\newcommand{\R}{\mathbb R}
\newcommand{\N}{\mathbb N}

\newcommand {\1}{{\bf 1}}

\newcommand{\F}{{\cal F}}

\newcommand{\s}{\sigma} 
\newcommand{\E}{\rm E}

\begin{document}
\begin{frontmatter}

\title{Mathematical Intuition, Deep Learning, and Robbins' Problem}

\begin{aug}
\author{F. Thomas Bruss}
\smallskip

\author{Universit\'e Libre de Bruxelles}
\end{aug}

\end{frontmatter}
 
{\bf Abstract.}
The present article is an essay about mathematical intuition and Artificial intelligence (A.I.), followed by a guided excursion to a well-known open problem. It has two objectives. The first is to reconcile the way of thinking of a computer program as a sequence of mathematically defined instructions with what we face nowadays with newer developments. The second and major goal is to guide interested readers through the probabilistic intuition behind Robbins' problem and to show why A.I., and in particular Deep Learning, may contribute an essential part in its solution.

This article contains no new mathematical results, and no implementation of deep learning either. Nevertheless,  we hope to find through its semi-historic narrative style, with well-known examples and an easily accessible terminology, the interest of mathematicians of different inclinations. 



\smallskip
\noindent
{\bf Math. subject index} Primary 60G40; secondary
68T01, 68T20.

\section{Introduction}

Many of us  remember that, when we were students, we heard ``Computers are important, but not for Mathematics''. These may have been the first instances (U. Saarbr\"ucken, around 1970) when I wondered to what extent this is true. Like many of us I felt that, in its essence, this should be correct.  Given the many new things of which we hear now
in terms of  A.I.-based approaches, would  we revise our point of view today?

\smallskip
The second time when this question came to my mind, and more forcefully,  was in  autumn 1976. We were chatting in the tea-room of the Statistical  Laboratory (U. Cambridge) when one of our lecturers came in and broke the news: ``The four-colour problem has been solved!" There was a joyful ``Oooohh !'' of true surprise and excitement in the room. It turned into a more modest Oo-oooo... when he added, after a short pause, ``... by a computer''.  Nobody said ``what a shame'', but a few of us seemed disappointed. 

\smallskip
Would this still be the case today?

\smallskip
Probably not. First, we are glad that the four-colours conjecture is settled.
Moreover, we agree that the proof of the completeness of  checking {\it all} cases  by a program was an achievement (Appel and Haken (1977)). No longer would it be seen as a  stigma that such a proof was different of what we were used to. Second, things would  be  likely to be different today since, over the years, we have learned. We always knew that we cannot compete with computers with respect to speed of computation. 
But now, we wonder in how far it is it just {\it this very speed} which  will break bounds in so many more directions than we can foresee.

To insist, let us imagine problems for which we see bounds which will prevent us to find a solution, such as intricate combinatorial problems, or analytical problems which suffer from the curse of dimension. But then, for which ones would we be sure that an enormous  increase of the speed of computation would not  break the bounds we see now? If somebody challenged us to give reasons why seemingly unsolvable problems will stay unsolved,  we would probably fail.

\subsection{The Namur experience}\label{NE}

The third time that I wondered about the importance of computers for Mathematics was some seven years later at U. Namur. Encouraged by theoretical results on optimal decision strategies ($1/e$-law, B. (1984), and  B. and Samuels (1987)), I wanted to test their interest for applications by computer simulation. Usually I do not do this, but I happened to have two interested diploma-students who were also very keen on programming, Rutten (1987) and Simonis  (1987).  
We wrote programs which would select (online) decision items from an incoming stream of an unknown number $N$ of rankable items arriving at (pseudo)-random times. Here is the essence: \subsection {A sequential decision problem}
Let $a, b\in R$ with $a<b$ real constants, and let $n \in \N$. Suppose $n$ points are distributed independently according to a continuous distribution function $F$ defined on an interval $[a,b].$ These points are ranked in some way. 
If $X_k$ denotes the $k$th observation then we call
$r_k$ its relative rank, i.e. its rank among $\{X_1, X_2,\cdots, X_k\}$. The objective of us players in the simulated game was to get online with maximum probability either a specific desired rank (such as for instance rank $1$) or, in extension, one in a certain range of ranks (as e.g. one of the
$q\%$ best ranks) under the rule that a point can only be accepted at the time of its arrival, and that exactly one of them must be accepted. What is the optimal strategy?

We competed with the computer in a series of games with the following set-up. The computer generated $n$ independent random variables according to the distribution function $F.$
Their increasing order statistics were interpreted as the arrival times $$T_1<T_2< \cdots <T_n $$  of the $X_j$ on $[a,b]$, all hidden to us at the beginning $a.$
At time $t$ the computer displayed \begin{align}\label{points}a~[.......\bullet.....\bullet...........\bullet...................\bullet..... |t ..........................]~b\end{align}
where the bullets stand for the successive arrival times on the time axis. 
Then the program  generated for each arrival time $T_j, ~1 \le j \le n,$ a uniformly distributed value  $X_j$ (not shown in the figure (\ref{points})). These values $X_j$ were ranked by the computer giving  to the smallest $X_j$ rank 1, to the second smallest rank $2$, etc.\footnote{One often ranks in opposite order, namely the largest value is rank $1$, the second largest value rank $2$ etc. Our ranking has here the advantage that ranks increase as the  $X_j$ increase. }    All what the screen showed to us at time $t$ was the  relative ranks of arrivals up to time $t$ with respect to preceding $X_j$'s. 

Thus, at time $a$ we saw in  (\ref{points}) no points at all, and at time $b$ we saw them all. 
 If the first four successive randomised values were for instance
\begin{align}\label{values}0.395..., 0.207..., 0.674..., 0.358...\end{align} then the above design would look in terms of relative ranks
  
 \begin{align}\label{Raenge}
 a~[.......1.......1...........3.....................2........|...........................]~b,\end{align}
 
 \smallskip\noindent
 since in (\ref{values}) the first is the smallest so far (by definition), so is the second, the third one is the third smallest so far, etc.  The only picture we saw at time $t$ was (\ref{Raenge}). The question was which relative rank we would select after time $t$ (say) if we can select at most one, and only directly {\it when it appears}. The idea behind the relative ranking was to mimic real-world conditions where one has no precise values attached to decision items but where one can still decide 'better-or-worse'.
  
 \subsection {Brain of the program}
 To give a specific example, suppose we have fixed at the beginning the objective "Get absolute rank 1." At the end (time $b$) we can see whether we
 fulfilled the objective. How should we play?
Note that a  decision problem of this kind captures a good deal of what happens in real-life. Opportunities pass and we would like to make an excellent decision, but we cannot return to opportunities passed over before.

For fixed $n$ the ``brain" of the program yielding the optimal strategy to get rank $1$  is simple, namely it is the one for a well-known best-choice choice problem (secretary problem). The optimal strategy to obtain overall rank $1$ is to stop on the  arrival time $T_k$ and select $X_k$ if
\begin{align}\label{SP} X_k {\rm ~has~relative ~rank~ 1, and ~} \sum_{j=k}^n \frac{1}{j}\le 1. \end{align} If there is none such $T_k$, the optimal strategy fails by definition.
\footnote{For other interesting objective functions  see e.g. Vanderbei (1980), Szajowski (2007, 2009), Bay\'on et al. (2018), Ramsey (2016). For modifications into several directions we refer to Goldenshluger et al. (2020).}

\subsubsection{Threshold strategies}\label{Thresh} A strategy of the form as in (\ref{SP}) is called a threshold strategy. $X_k$ is the first observation (if any) which has a certain property.  In (\ref{SP}) the required property is that $X_k$ satisfies a quality constraint (here meaning that is has relative rank $1$)
and that the index $k$ (which measures here the arrival number)  satisfies a time constraint, namely the given sum constraint in (\ref{SP}).

 It will be convenient throughout this paper to call strategies to stop the first time that an observation satisfies a given quality constraint,  also in continuous time,  {\it threshold strategies}.
\subsubsection{More general objective functions}
For the more general objective of getting one of the  $q$ percent best absolute ranks for some $q$), a nearly optimal strategy is also still easy to program. Moreover, recall that a continuous distribution function $F$ sends all arrival times into uniform random variables on $[0,1]$, that is, \begin{align} X_1, X_2, \cdots X_n {\rm ~i.i.d. ~}F{\rm-distributed~and~ }F{\rm ~ continuous}~~~~~~~ \nonumber\end{align}{\rm~then~}
\begin{align}\label{Transf} F(X_1), F(X_2) \cdots , F(X_n) {\rm~i.i.d. ~uniform  ~on} ~[0,1].\end{align}
 
Hence the $F(X_j)$ are all i.i.d  $U[0,1]$ random variables, and all what we need for applying a desired threshold strategyis the table of threshold arrival times in $[0,1]$-time, a memory stack for relative ranks of arriving options, and a simple integral operator applied to the observed empirical arrival time distribution. (See B.(1992), p. 13, theorem 5.1, and p. 33, Table 1.) Hence, clearly, whenever we speak of a continuous arrival time distribution $F$ on some interval $[a,b]$ we can assume w.l.o.g. that it is the uniform distribution on $[0,1]$.

 (For a different approach see e.g. Villeneuve S. (2007).)

\subsection{Statistical inference and the power of speed} \label{SIPS} To increase the appeal for applications, the number of observations  $N$ can  be randomised in each run. Unless $N$ turns out very small, the optimal threshold times for relative ranks hardly change  as the theory shows nicely (B. 1984, equation (3)). In a further generalisation, if the computer randomises the arrival time distribution function $F$ in some {\it finite}  basket ${F_1, F_2, \cdots , F_k}$ of such (continuous) distribution functions, then, if we knew which one the program  took, we could use the same table of threshold times for ranks in $[0,1]-$time. 

Statistical inference comes in with the  idea  to act as if the chosen distribution function were, at time $T_j,$ the function $\tilde F_j$ which is, in some sense, closest to the observed {\it empirical} distribution function $F_{T_j}^{\rm emp}(t)$ of arrival times on $[0,T_j].$ If I remember well, I suggested to use integrated squared distance to measure closeness, i.e. to put
at the arrival time $T_j$ \begin{align}\label{distemp}
\tilde F_j(t) :=\arg \min_{G\in \{F_1, \cdots, F_k\}} 
\left\{\int_0^{T_j}\left(G(u)-F_{T_j}^{\rm emp}(u)
\right)^2du \right\}, j=1, 2, \cdots \
\end{align} and to set, for $t\ge T_j$,  $F(t):= \tilde F_j(t)$. Thus the tail of $F$ may change with each new arrival. This is of course no problem  for a computer, but it is so for a human being trying to keep up.

In a fair competition, both the computer and the player have to guess online the chosen distribution function and then also to deduce accordingly information about the likely range of $N.$ We randomised $N$ (uniformly) in $\{1,2, \cdots, M\}$ and applied the correspondingly optimal strategy of choice for the fixed objective. 

Unless we confined to small $M$, the computer would beat us in playing. No chance to compete with the speed of the computer to concentrate on "the most likely" chosen distribution function and to estimate $N$ accordingly. This was fully compatible with our intuition and proved indeed to be the computer's dominant advantage.

\section{Steps into A.I.}\label{compatibility}\label{SBE}
Finally, some two years later (U. of Cal. Santa Barbara) I picked up, just for fun, the programs with a different objective. My wish was  now to modify them such that the computer would play well against me even if I kept my objective for the games secret (!)

The simple trick was to let the ``uninformed'' program use incidence matrices of compatibility which kept track of my decisions. 
I forgot the details, but here is the philosophy. If I made a choice the program would "guess" that (true or not) this choice was not too far from my secret goal. As explained before, a single sufficiently fine table of optimal thresholds (in uniform $[0,1]$-time) was sufficient to deal with objectives of getting a specific , or alternatively, a rank in the
set of the $q$ percent top ranks.  

For example, if the player's choice in a game turned out rank 3, then this was compatible
with the goal being to get one of the three best, i.e. rank 1 or 2 or 3, but also compatible with the goal of getting one of the 10 percent best if $N$ turned out at least $30.$ Of course, this was not conclusive at all. I may have intended to get the rank 1, but was just unlucky, or, on the contrary, I may have planned to get one of the best 20 percent and did, just by chance, well. However, this was enough to give the computer information in form of neighbourhoods of my likely objective. The constraint was of course that I not was allowed to change in the {\it same} series of trials the objective, or the basket of possible distribution functions.

\subsection{Goals, neighbourhoods and convergence}
If one reads about deep learning, one may see remarks of the form "We understand why our method works, but we do not know why it works so well". Mathematicians may object that if one {\it really} understands  the "why" then one can also understand  the "why so well". 

In the above described case, learning specific goals meant learning to narrow down neighbourhoods of them, and in a sufficiently long series the program would learn a reasonable neighbourhood of the likely goal quickly enough. Since the objective functions are rather robust, and since the computer  assessed the most likely $F_j$ in my (finite) basket by computing $\tilde F$ above
and reading off the optimal threshold times so much quicker than I could do, the program often beat me in the longer run. Only if $N$ turned out  small
I had an advantage by knowing what exactly my goal was. However, if $N$ is uniformly randomized on $\{1, 2, \cdots, M\},$ then $N$ is
rarely small if $M$ is not small. 

Thus robustness and the speed of convergence
to a narrow set of goals explain why our program did "so well".  A good deal of this robustness is pointed out in B. and Samuels (1990). (For more details on so-called fine sets for vaguely defined objectives, the interested reader is referred to B. (1992), pages 23-28.) 

In summary, I was the father of  the decision program, but I had no chance to beat my child, even though it did not know what my goal was.
My feeling at that time was that I had made a step into what I would call now  A.I.  (not yet  deep learning.) It  was based on a simple idea, but the outcome was impressive. 
My conclusion was
that in the future I should be less impressed by A.I. than by achievements of  hardware and software engineers
who have made this possible.  
 
 \smallskip I wonder now, whether this conclusion was again premature.
Will the new things of which we hear in A.I. trigger off for us a new way of thinking?  What I call ``new'' is the surprising developments in machine learning, and specifically in  {\it deep learning}. {\it New} lies in the eyes of the beholder, of course, and the author, likely as many other readers, does not know enough about deep learning.
The type of A.I. we all understand without further preparation is what we can write down in a flow-diagram style, where a {\it condition} is a mathematical statement, where a {\it  Goto} is a well-defined instruction, and where an iteration has a definition of where it begins, and where it ends. 

Yes, deep learners should pardon us, but the offspring of the Algol 60/68 generation still do exist.

\section{A.I. Definitions and common terminology} 

Readers with some experience in A.I. can skip this section and pass to Section 4.
For all others we review briefly important definitions and common terminology in A.I. we use and exemplify them, when possible, by what we described in the examples treated in our Introduction. (Guidelines for further reading may be found  e.g. in Goodfellow and al. (2016), and Acemoglu and Restrepo (2019))

\smallskip
{\it Artificial intelligence} is the generic name  for the large domain of methods and tools to provide intelligence to machines  (usually computers) which can mimic or outdo human intelligence. Some authors using the word A.I. in this domain would include the study of relevant hardware, others would not. 

There is no clear 
agreement where A.I. begins. The solution of the four-colour problem cited in Section 1 for instance is often seen as an early convincing example of A.I. However, we recall that the program is a skilful well-organized  scheme of exhaustion to check all possible cases. The fact that certain parts of this proof can be made more elegant by applying A.I.-tools is something different. Thus the four-colour problem need not be seen as an early step into A.I.

\smallskip
{\it Turing test}~ The Turing test, named after Alan Turing  (1912-1954), is an intuitively appealing  test to decide whether the type of intelligence we encounter in a specific case should be considered as A.I.  In simplified form it says:
If a computer and a human player, say, interact, and a neutral observer,  who can only see the interactions of the two, cannot decide who is the human agent, then the computer demonstrates A.I.. The Namur-experiments (see \ref{NE}) (fixed objective) would pass this test. The first version would not  pass it when the human agent is allowed to keep the objective secret. With the additional  modification based on keeping track of the player's decisions (compatibility matrices) it would pass it in most cases.

\smallskip
{\it Machine learning}~~The word {\it learning} is understood as it is in everyday language. Machine learning says not more than that this is done by a machine. A finer distinction comes with the task to be executed by the machine, and with the tools it can use. For instance in the first part of (\ref{NE} ) we would not speak of learning. It suffices to know the input $n$ (the total number of variables) and then to apply
the deterministic algorithm (\ref{SP}) to obtain the optimal solution.
The same would be true if the machine should solve similar problems for more general objectives. Only the corresponding algorithms would be more complicated. (Approximate optimal stopping was studied in papers by K\"uhne and R\"uschendorf (2000, 2003). The connection with optimal stopping problems in some more generality is treated in Christensen and al. (2013), Villeneuve (2016), R\"uschendorf (2016), and Li and Lee (2023).)

The need of learning is evident, if the number of variables becomes a random variable ($N$) and the more so, if the distribution of arrival times $F$ is itself a random variable with values in some finite set
$\{F_1, F_2, \cdots, F_k\}$. The machine is then bound to have to learn from what it sees. We suggested sequential estimation of the true $F$ by minimizing sequentially (\ref{distemp}). In the present paper we gave no further details, but it is then a classical step in Statistics to infer, again sequentially, on the unknown $N.$
The difference compared with the problem for fixed $n$ is that now, to apply the procedure, the machine must learn the data (arrival times), and  infer from the data. 

\smallskip

{\it Supervision in machine learning}  Machine learning may display several levels of supervision. For example, in our setting of machine learning (with unknown $N$ and unknown
$F\in \{F_1, F_2, \cdots, F_k\}$) we may see the functions $F_j$ also as data  (in a wider sense). A.I.-specialists may then speak of {\it supervised} learning, since the machine
uses {\it labeled data} (the  $F_1, F_2, \cdots, F_k$) 
together with the data (in a strong sense) which are the arrival times in each run. It is useful to keep
the meaning of {\it  labeled} in mind as data to which the learning computer is directed in some sense.

If only partial information about labeled data is given  (for example $F$ lies in some subset of $\{F_1, F_2, \cdots, F_k\}$) one may 
speak of a form of {\it semi-supervised} learning, but clearly many more different cases can then be distinguished. 

\smallskip
{\it Deep learning}   Deep learning is a subset of machine learning that focuses on so-called {\it neural networks.}  Neural networks are
seen as the {\it black boxes} of a deep learning program. The word deep refers to the number of {\it inner layers} in network. The higher this number of inner layers, the larger can be the complexity of the problem, and usually one speaks of "deep" if it is at least two. Here the two external layers input and output are not counted.  A beginner's guide can be found in Nicholson (2020). See also Wang et al. (2020). Deep learning problem specifically for optimal stopping problems were studied in Fathan A. and Delage E. (2021).

 \smallskip
 {\it V-learning} In V-learning (for a recent contribution see Li and Lee (2023)), the goal is to estimate the value function without explicitly learning a policy. Instead, it focuses on learning the values of states and state-action pairs directly.
The V-learning algorithm is based on the Bellman equation  and reinforcement learning. The Bellman equation expresses the value of a state as the expected sum of rewards that can be obtained starting from that state and following a given policy. 

\smallskip

{\it Deep reinforcement learning}
 In none of our examples given in our Introduction we worked with neural networks, and the author has no experience to offer in deep learning. However, deep learning and reinforcement will play later on an important role in what we will suggest as a way to attack Robbins' problem. There we will try to show that intuitions can help us to funnel the attack through deep learning in a promising direction.

\section{Deep problems}

The idea of deep learning is that the machine would learn from data. If we are allowed to understand the set-up of learning as a set of  programs which define Goto-instructions as functions of data, and run iterations on this set of data, possibly augmented by newly obtained data, etc,  then, as we said already, we understand. Deep learning however is something different, and for this the set of data enabling learning  (usually called {\it training set}) should be large. 

Actually, deep learning techniques require huge training sets of data. This seems intuitively clear. 
But, then,  a caveat. In several cases, we do not see where these suitable training sets should come from.
To discuss this,
we exemplify a few famous open problems.

\subsection{Number Theory} Let us take the Twin-primes conjecture as our first example. We have known since the times of Euclid, that there are infinitely many primes, and we still love Euclid's short and elegant proof.  But then, what would be a training set for the twin primes other than a set of positive  integers spanned by infinitely many "candidates"?

Also, is the intrinsic problem of imagining training sets not the same in other examples of celebrated open problems, such as e.g. the Goldbach conjecture, the Riemann hypothesis, and so many others? And then, for instance, what could be training set for the Collatz conjecture ($3n+1$ problem) other than the set of all positive integers?

In the case of the Twin-primes conjecture we have intuitive (pseudo-probabilistic) arguments based on the Prime Number Theorem (e.g. De la Vall\'ee Poussin (1897)), why it should be true.  
The same is true for the (strong) Goldbach conjecture, and this even in a reinforcing sense. If it is true for $n$
then it is (a prior, not talking about pair/impair) probably true for $m$ with $m>n$ since the number of prime-candidates $b$ and $m-b$ yielding the sum $m$ is  increasing. 

A very crude and well-known version of the heuristic probabilistic argument for the strong form of the Goldbach conjecture is as follows. According to the Prime Number Theorem, an integer $m$ selected at random has roughly a  chance of $1/\log(m)$ being prime. If $n$ is a large even integer and $m$ is a number between $3$ and $n/2$, then, under the assumption of independence, the probability of $m$ and $n-m$ simultaneously being prime is in the order of $(\log(m)\log(n-m))^{-1}.$ Hence one may expect the total number of ways to write a large even integer $n$ as the sum of two primes to be roughly
   \begin{align} \sum_{m = 3}^{n/2} \left(\frac{1}{\log(m) \log (n-m )}\right) \to \infty~{\rm as ~}n\to \infty.
 \end{align} Here we neglected interdependencies between occurrences 
 of primes, but
this would not question the general reinforcement effect of growing even numbers.

Pseudo-probabilistic proofs are no rigorous probabilistic proofs, such that intuitions need not be  helpful to invent adequate training sets. \footnote{However, some models explain very well why certain results in Mathematics must hold in great generality. See e.g. the elementary approach to Taylor's polynomial in B. (1982).}
Moreover, we probably would agree that, if we were able to construct a suitable training set for a given problem, then we would  have  enough insight to make progress on it anyway. Hence our  conclusion: Deep learning is unlikely to be helpful in these famous problems. 

Are there common features in famous open problems?

\smallskip
Some are evident. In our examples, the answer has to be of the form true or false, that is "Yes or No". There is nothing in between Yes or No. 
Also, many hard problems in number theory contain a  $\forall$-statement. In the Twin-primes example, we would like to know, true or false,
\begin{align} \forall n \in \N ~~\exists p\in \N, \,p>n :~ p{\rm~ and~} p+2 {\rm~are~prime~numbers.}\end{align} The implication of this is also a common feature. If suitable training sets did exist, they would have to be arbitrarily large subsets of the positive integers.

Consequently, yes-or-no questions involving {\it for all} statements in some form or other in open problems are probably no promising candidates for profiting from learning algorithms. (Since we cannot have a both solid and useful definition of an "open" problem there exists probably no  useful theorem which puts this statement on a solid ground.) 

In agreement with the title of the present article, we now take the liberty to return to intuition and look at examples of deep open problems, where the restrictions imposed by {\it for all statements} are less severe. In this vein  we leave Number Theory and turn to a  problem which starts like a problem in elementary analysis. 

\section{Analysis}\label{ElAnal}
 Let $(\ell_n)_{n=1, 2, \cdots}$ be a sequence of real numbers satisfying

\smallskip
(i) $\forall n\in \N: \ell_n \le \ell_{n+1}$

\smallskip
(ii) There exist bounds $L,  U \in \R$ such that $L\le \ell_n\le U$ for all $n\in \N.$

\smallskip\smallskip
\noindent {\bf Questions}: 

\smallskip
Q1 ~Does $\ell=\lim_{n\to\infty}\ell_n$ exist? 

\smallskip
Q2  ~If $\ell$ exists,  what is its value?

 \medskip
 \noindent Q1 is trivial  (Bolzano-Weierstrass), and Q2 (a priori) meaningless, of course. We see from (i) and (ii) that $(\ell_n)$ is both increasing and bounded above, and thus $\ell$ exists and must satisfy $\ell\le U.$
The lower bound $L$ plays here no role at all.

In order to answer Q2
we would need additional information on the  $\ell_n.$  
Let us recall a few cases where things work out.
If, for instance,  the $\ell_n$ are defined recursively by
$\ell_{n+1}= F(\ell_n)$ for some continuous function $F:\R \to \R$ then we know that we can obtain the limit $\ell$ if the equation $F(x)=x$ allows for a solution $x\in[L,U].$ The same stays true for subsets of $\R^k$ if we have a recursive definition of the form \begin{align}\label{Fkl}  \ell_{n+1}=F_k(\ell_n, \ell_{n-1},\cdots , \ell_0):= F_k(\ell_n, \ell_{n-1},\cdots, \ell_{n-k+1})\end{align} for some fixed $k\in \N$ and some continuous function $F_k:\R^k\to \R$ because the uniqueness of the limit, if it exists, implies that for a fixed $k\in \N$, the limit $\ell$ must solve \begin{align}\ell=F_k(\underbrace{\ell, \ell, \cdots, \ell}_k).\end{align}

We can go on  with similar fixed-point arguments if we allow the recursive definition to depend on {\it all} preceding terms by requiring in addition that the influence of the ``earlier" part of the history in \ref{Fkl} {\it fades away sufficiently quickly} as $n\to \infty.$ 
However, with this,
we have arrived more or less at the end of getting an answer from a fixed-point equation. And this is the point. We are 
lost if we cannot describe, at least to some extent, the dependence of the variables.

We know many domains in Mathematics where we can find such examples. We choose a specific problem
which is on the borderline between analysis and probability theory:

\section{Robbins' problem}
Our choice is a celebrated open problem in the domain of Optimal Stopping. It is known as Robbins' problem of minimizing the expected rank. 
(This is, by the way, the same Herbert Robbins who co-authored with Courant {\it What is Mathematics?}' (Courant and Robbins (1941)).  

 \smallskip
 Robbins  presented this problem at the end of his memorable talk on the International Conference on Optimal Stopping and Selection in 1990 (U. of Massachusetts). This problem is less known than the famous problems we mentioned before. However,  it  is deep and is an easy-to-decribe representative in a whole class of problems in Probability, about which we do not know what to do. 
\subsection{Definition of Robbins' Problem}
Let $n\in \N$ be a fixed, and let $X_1, X_2, \cdots, X_n$ be independent and identically distributed random variables uniform on $[0,1].$ We can observe them sequentially, i.e. in the order $X_1, X_2, \cdots,$ and we must select {\it exactly one} of them. A selection is only possible at the time of observation, and if $X_k$ is selected, then the  decision is irrevocable
and the game is finished. At step $n$ we see the whole picture, and if we have selected $X_k$ then we occur the  loss
\begin{align} L_k=\sum_{j=1}^n \1\{X_j \le X_k\}.\end{align} In other words, if we denote the increasing order statistics of the $X_j$ by
\begin{align}X_{1,n} \le X_{2,n} \le  \cdots \le X_{n,n}\end{align} and if we have chosen $X=X_{k,n},$ then our loss is $k.$ Hence, here is the problem different form the best-choice problem we mentioned before. We see the {\it values} and not only their relative ranks, and the {\it loss is the rank} of the accepted observation !  At time $n$, the random variable $L_k$ becomes deterministic, taking the value $R_k$ which is the final rank of $X_k$
among the whole sample $X_1, X_2, \cdots, X_n$. If we must choose exactly one variable,
what sequential strategy will minimise the expected loss? \footnote{Robbins announced it by saying ``Finally, here is the problem which I'd like to see solved before I die" and he looked into the audience in a way we all felt he meant it. So much for the name. Sadly, Robbins'  wish did not realize. He died February 12th, 2001.}

\section {A brief analysis of Robbins' problem}\label{BAna-RP}
Since for each $n$, a strategy results in at most $2^{n-1}$ yes-or-no decisions (select or go on), there are only finitely many relevant strategies. Thus for each $n,$ an optimal strategy must exist. Let $v_n$ denote the corresponding optimal value, i.e. the minimal expected (final) rank.

\smallskip
The following results are proved in B. and Ferguson (1993) 

\smallskip
(i) $v_n$ is increasing in $n.$

\smallskip
(ii) The sequence $(v_n)$ is bounded below (trivially by $L=1$)
and bounded above by another known value from a related problem studied in Chow et al. (1964), which is $3.869\cdots$. Hence its limit \begin{align}v=\lim_{n\to \infty} v_n\end{align}exists.
Thus so far, the sequence $(v_n)$ satisfies the setting for the sequence $(\ell_n)$ defined in Section (\ref{ElAnal}). Unfortunately we do not see an easy relationship between the $v_n$. Additional information is  harder to get, but the next one 
is both simple and important.

\smallskip
 (iii) Smaller $X_k$ have smaller ranks. Hence, it is intuitive that the values $X_k$ and the corresponding final ranks $R_k$ (which will be known at time $n$ only) should be positively correlated. Recall that {\it correlation} is a measure of dependence of one random variable of another one. We can compute it  and obtain (B. and Ferguson (1993), (1.6) - (1.8))
\begin{align}\label{corr}\forall \,n\in \N,~ \forall \,1\le k\le n:  {\rm corr}(X_k, R_k)= \sqrt{\frac{n-1}{n+1}} \to 1 ~{\rm as ~}n\to\infty.\end{align}

\smallskip
(iv) The strong positive correlation in (\ref{corr}) suggests to propose a strategy
by just looking at time $k$ at the observed value $X_k$ and to select it if and only if $X_k$ is smaller or equal to some threshold $\varphi(X_1, X_2, \cdots, X_{k-1}; X_k; n)$. This is a threshold strategy. If moreover, we ignore at each step $k$ all preceding values $X_1, \cdots, X_{k-1},$ then we speak of a {\it memoryless} threshold strategy, or in short  {\it ml}-strategy. We denote the optimal value obtainable for $n$ observations in this restrained class by $\tilde v_n.$ (We will return to
$ml$-strategies in more detail when needed later on.)

\smallskip
(v) Seeing the numerical values of the random variables allows us to rank them. Numerical values provide thus at least as much information as rank-information. Hence one cannot do better with rank information than with seeing the $X_j.$,  We can similarly show that this bears over to memoryless strategies, so that the sequence $(\tilde v_n)$ is increasing and bounded above. Hence \begin{align} \tilde v=\lim_{n\to \infty} \tilde v_n~{\rm exists~ and~} v=\lim_{n\to \infty}v_n \le \tilde v.\end{align}

\smallskip
(vi) Lower bounds of different levels for $v_n,$ and thus for $v$ can be obtained by a truncation argument (B. and Ferguson (1993)). We say we truncate the loss at level $j$ with $1 \le j\le n,~ j\in \N,$ if the loss generated by $X_k$ is defined to be equal to $\min\{j, L_k\}$. Therefore, clearly,
\begin{align}\label{BAna-RP(vi)}\forall 1\le j\le n: v_n(j) \le v_n.\end{align}
Truncating at the level $j=5,$ B. and Ferguson obtained the lower bound $L\approx 1.908$. This lower bound can be improved slightly but not by much (see (B) below).
\section{The intrinsic difficulty of Robbins' problem.}
To find $v_n$ and $v=\lim_{n\to \infty}v_n$ we started with a clean analytical approach. So then, why do we  stop here?

\smallskip
There are two reasons. 

\smallskip
(A) The first reason is that we seemingly cannot compute
$v_n$ for large $n.$ We do not know the optimal strategy for larger $n$. It is trivial for $n=1$, almost trivial for $n=2,$ and still easy for $n=3$, but that's it as far as we can say "easy".  The problem is that the optimal strategy is {\it fully history dependent}, which means that the optimal strategy depends at each step $k$ on the complete preceding cloud of values  $X_1, X_2, \cdots, X_{k-1}$. (The  order in which the points arrived is of course irrelevant.) 
Since the proof of full history-dependence is somewhat involved, we refer for it to B. and Ferguson (1996), pages 9-13.
  
  \smallskip
From a decision-theoretical point of view, the interpretation of full-history dependence is that there is no sufficient statistics for optimal decisions other than the whole history (cloud) itself.

To understand the difficulty to compute $v_n$ precisely, we refer to the figures Figure 1-Figure 4 of Dendievel and Swan (2016). These authors coped with the challenge to solve the problem for $n=4$. Their graphs for the composed acceptance regions of the corresponding optimal strategy are complicated and show so little structure that one sees little hope to compute $v_n$ precisely for $n\ge 5.$ 

 \medskip
(B) The second reason is equally bad news. The idea
to compute for, larger $n$, the $v_n$  by the truncation method described in (iv) of Section \ref{BAna-RP}, which comes naturally to our mind, is just hopeless. As shown in B. and Ferguson (1993), the storage demand to implement the truncation method increases exponentially in both the number $n$ and the truncation level $j.$ Under the hypothesis that the capacity of computers increases exponentially, this is far away from the double-exponential growth we would need.

\smallskip
Summarizing the impact of (A) and (B), many authors have stopped working on Robbins problem when realizing how much damage the intersection of both does to any hope to solve the problem. Indeed, not much was contributed after the nineties. Even the following question is still open:

\section{Do $v$ and $\tilde v$ coincide?}
This is a yes-or-no question. However, since we only speak about the limits, there is {\it no for all}-statement in it. Finding one upper bound $U$ for $v$ with $U<\tilde v$ would do. Hence the main open part of Robbins' problem is submitted to a weaker constraint than all the celebrated problems we mentioned before, and thus, at least conceptionally, easier to answer.

If $v=\tilde v$, then  $v=\tilde v\approx 2.327$ (Assaf and Samuel-Cahn 1996.)
However, why does the question "$v=\tilde v$?" given what we know about Robbins' problem, becomes an important one?
We recall that $\tilde v$ is the value of the limiting optimal {\it ml}-strategy as $n\to \infty.$ We also remebr that the best known lower bound for $v$ is $L=1,908...$

It may come as a surprise that, with a gap of  $G=U-L\approx 0.419 $ one cannot say more.  We can interpret this by saying that $(U-L)/(L-1)>0.45$, so that at least $45$ percent of what may be  achieved is seemingly lost due to insufficient intuition.
Interestingly, without going here into details, we can show that for all $n$, one can find strategies which do strictly better than the optimal {\it ml}($n$)-strategy.  However, we do not know whether the difference will {\it persist as $n\to \infty.$} 
Note the vicious circle. We know that the limit $v$ exists and that it lies somewhere between $L=1.908\cdots $ and $U=2.327\cdots.$ However, for larger $n$ we cannot compute the optimal strategy and its value $v_n,$ and hence we cannot decipher how the difference $\tilde v_n-v_n$ behaves. 
As so often in analysis, one would like to know then at least its limiting optimal value $v.$

\begin{figure}[ht]
\centering
\includegraphics[width=0.60\textwidth, angle=0]{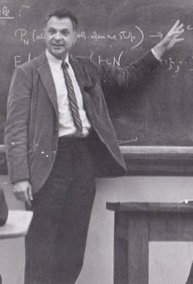}
\centerline {Figure 1}
\begin{quote}
Photo of Herbert Ellis Robbins taken in 1966 during a talk at Purdue University.   \end{quote}
\end{figure}

\subsection{The importance of knowing $v<\tilde v$ or $v=\tilde v.$}

For a probabilist, knowing whether $v <\tilde v$ or $v=\tilde v$ is almost as enticing as knowing the precise value of $v$. 
How come? 

It is the same  as what is behind the feeling of almost every mathematician. If we have a proven result  (we have not)
then it is the joy of having proved a mathematical result. Otherwise, we would appreciate a true understanding of how far one could possibly go, and what our intuition could contribute to this understanding.

To be explicit on  this with respect to Robbins' problem:  We understand where $\tilde v$ is coming from and what it means , namely it is the best limiting value which we can obtain in the class of $ml$-strategies.  If $v=\tilde v$ then this means that full-history-dependence is, as $n\to \infty$, of no importance at all, ands it would be neither for possible applications.  The strong correlation of an observation $X_k$ with its final rank $R_k$ combined  with the strong law of large number of the $X_k$  are together enough to imply that $v$ and $\tilde v$ coincide.

But what if $v< \tilde v$? Indeed, many people who have worked on this problem believe that this is the case.  All we know for sure so far is that $1.908 < v \le \tilde v \approx 2.32$, and the lower bound is hardly of interest. The point is that if $v < \tilde v,$ then  we do not understand the problem as well as we should, because, if we did, then we would know how to improve on {\it ml}-strategies  such that the improvement would not fade away
 as $n\to\infty.$ Where does our intuition fail to grasp essentials?
Several people have worked on this.
Published contributions were made by Assaf and Samuel-Cahn (1996), Gnedin (2007), 
Gnedin and Miretskiy (2007), Meier and S\"ogner (2017),  B. and Ferguson (1993, 1996), B. (2005), and others. Robbins himself had said in his talk that $v$ is ``not so far from 2",  but we never could find out what he meant by `not so far'. \footnote{According to S.M. Samuels, L. Shepp, and D. Siegmund, (private communications), Robbins was usually generous in sharing scientific results, but on "his" problem he seemed particular. }

\section{Logic and Intuition}
Suppose that $v<\tilde v$. Then the limiting optimal strategy is clearly different from the limiting memoryless threshold strategy. Hence there must be cases where optimal behaviour should overrule the instruction given by the optimal $ml$-strategy.  What exactly is it that our probabilistic intuition does not to grasp?

If $v<\tilde v$ then the answer must then be hidden in one of two possible scenarios. At certain arrival times $t$ the current information (i.e. the "cloud" of points in $[0,1]$ which have appeared already) should persuade us not to stop and go on, although the current observation is {\it smaller} than the optimal {\it ml(n)}-threshold, or else it should persuade us to stop although the current observation is {\it larger} then the optimal {\it ml(n)}-threshold.
In a colloquial terminology we may speak of dissuading ``pre-clouds" and persuading ``post-clouds" for stopping. 

\smallskip
Here is an illustration. (Clearly, for large $n$, only neighbourhoods of $0$ are of interest for most observations.)

\begin{align}\label{prep1}~~~~0~[ ~ .....x.......x.x.xx.\bullet ..t_n^*...|.....|.......|...|...........|...|.............\end{align}
\begin{quote}
{\it Dissuading pre-cloud} ~The $\bullet$ denotes an observation at time $t,$ which is smaller than the optimal $ml(n)$-threshold $t_n^*$.  In the above figure, 5 values (denoted by $x$) precede it, and none of these $x$ was accepted when it appeared. By accepting  the $\bullet$ the current pre-cloud would push its rank to at least 6. It may be better to go on.
\end{quote}
\begin{align} \label{prep2}~~0~[ ~ ............. ,t_n^* . \bullet.x.x...x x.. .......x....|......|..........|..|...|.........\end{align}
\begin{quote}{\it Persuading post-cloud} ~In the same spirit as above, the $\bullet$ may now be a good candidate for acceptance, even though it is larger than the optimal $ml(n)$-threshold $t_n^*$, simply because it just escaped
the push-up of the rank by 5 close values (denoted by $x$)  thereafter. 
\end{quote}

But now the same kind of reasoning would again pose new questions of all kinds. So, for example, how many values in a post-cloud would compensate at time $k$ for a given pre-cloud in an interval of the same length? Should we also take combinations of the following form into account?

\begin{align} \label{prep3} ~~0~~~[ ~ ...{\bf[}...x..xxx........... ,t_n^* . \bullet..x.x.xx..x .......x...{\bf]}....................\end{align}

\begin{quote} If the $ [...\bullet]$ and the $[\bullet ...]$
in the preceding figure are supposed the be of the same length, is it then the difference of the $\#x$ to the left, respectively to the right, which should be the relevant parameter? But then, how many as a function of different lengths, etc? 
\end{quote}

Answering such questions by defining the parameters  will depend on what we call here {\it sequential calibration}.
This calibration allows for so many possibilities that we feel (actually for the first time in the present paper) a need
to have it done by an automatism. 

\smallskip
Deep learning specialist will have ideas which are probably independent
of ours, guided by intuition. But we equally respect intuition and take the liberty to make suggestions
in the spirit of reinforced learning.
	
\subsection{How could reinforcement learning work?}
The general idea of reinforcement learning is as follows:
Developers devise a method of rewarding desired behaviours and punishing negative behaviours. This method assigns positive values to the desired actions to encourage the agent and negative values to undesired behaviours. This programs the agent to seek long-term and maximum overall reward to achieve an optimal solution. So the desired result would be that a strategy of selection different from the optimal $ml(n)$-strategy leads for a large $n$ to an expected rank strictly below the upper bound $U.$
These long-term goals help prevent the agent from stalling on lesser goals. With time, the agent learns to avoid the negative and seek the positive. This learning method has been adopted in artificial intelligence as a way of directing unsupervised machine learning through rewards and penalties.

It becomes evident what we have in mind. A.I. should give instruction to overrule the optimal $ml(n)$-strategy by refusing
at time $k$ a value below the optimal threshold $t^*:=t^*(n,k)$ if dissuasion at time $k$ by a pre-cloud is sufficiently strong (and go on), and accepting a value above 
$t^*(n,k)$ if persuasion  by a post-cloud is sufficiently strong. At time $n$ (the end) the final rank will be noted. To enable this strategy we need to know at least a very good approximation the optimal $ml$-strategy and the expected value we can obtain with it.

\smallskip

\subsection {Approximation of the optimal $ml(n)$-strategy and calibration}

Let $\{\varphi[n]\}=\{\varphi_1, \varphi_2, \cdots, \varphi_n\}$ be the $ml$-strategy to accept the first $X_j$ with $X_j\le \varphi_j.$ We recall that we can confine our interest to i.i.d. uniform random $X_j$ on $[0,1],$ and thus also on
on threshold values $\varphi_j\le 1, j\in\{1, 2, \cdots , n\}$. To assure that one $X_j$ is accepted, we put $\varphi_n=1$.
Applying  this strategy we accept $X_j$ if $X_1>\varphi_1, X_2>\varphi_2,\cdots , X_{j-1}>\varphi_{j-1} {\rm ~ and~}X_j<\varphi_j$
which, due to the i.i.d. assumption on $[0,1]$,  occurs with probability
\begin{align} P(X_j~{\rm~is~accepted})=\varphi_j~ \prod_{\ell=1}^{j-1} (1-\varphi_\ell).\end{align} Also, if $X_j$ is accepted,
its expected final rank will be its relative rank $r_j$ plus the expected number of later observations smaller or equal $X_j$,
i.e. \begin{align} \E(R_j | X_j {~is~accepted}) = r_j + (n-j) X_j. \end{align} Hence we have  (in principle) all the ingredients to compute the precise formula of the expected value obtained under  $\{\varphi[n]\}$ and then to minimize with respect to $\{\varphi_1, \varphi_2, \cdots, \varphi_n\}.$ This is technically more involved, but the following approximation is expected to be good enough for our needs.
Putting for $c>1$ 
\begin{align}\label{approx}  \varphi_j (n) = \frac{c}{n-j +c}, ~ 1\le j\le n,\end{align} allows to compute  the corresponding expected rank easily. Minimizing the latter with respect to $c$ yields $c=1.9469...$ and  the corresponding expected rank becomes $2.3318 \cdots$. The best known proven interval for $\tilde v$ namely $2.295<\tilde v<2.327$ is due to Assaf and Samuel-Cahn (1996). Since we are mainly interested in $v$ (and not $\tilde v$) we want to stay compatible in the comparison criterion and thus propose

\begin {quote}\label{crit}
\smallskip
1. to use the $ml$-threshold functions $\varphi_j(n)$ defined in (\ref{approx}) with $c=1.9469,$ 

\smallskip
2.  to take (for very large $n$) the expected value $2.3318$ as $U$, i.e. as our point of reference.
\end{quote}

If the average of accepted ranks in a large number of simulations is above
$U$ we discourage whereas we encourage if it is below $U.$ This can be done for each run within a single set of runs. A strategy $S$ is seen as better than a strategy $S'$
if the average of the obtained  ranks in a run  is lower than the one for $S.$
of simulations. Here is a suggestion of relevant parameters. 

\begin {quote}\label{Chall1}
$n=$ number of observations (large); 

$k=$ index (time) of observation  $X_k; 1\le k\le n$;

\smallskip

By definition of $n$ (total number), and $k$ as the current observation number. Hence these two parameters are not affected by learning.

\smallskip

$d_{pc}=$ length of the observed dissuading pre-cloud interval $[X_k-d_{pc}, X_k]$; 

$N_{d_{pc}}(k)=$ number of observations therein at time $k$

\smallskip

$p_{pc}=$ length of the observed persuading post-cloud interval $[X_k,X_k+p_{pc}]$; 

$N_{p_{pc}}(k)=$ number of observations therein at time $k$

\smallskip
Variations of these parameters
\end{quote}

\smallskip
\noindent Here is a suggestion to funnel reinforcement learning through
a variation of the winner's rule, although the author would have no particular reason to defend his
choice. 

\medskip
\begin{quote}
	{\it Play a randomized winner's rule}
\smallskip

	(I) If the strategy pursued in the last run yielded a rank smaller then $U$, then use the same strategy for the next run.

\smallskip	
	(II) If not then randomize (with probability $1/2,$ say) to decide whether to repeat the strategy applied in 
	last run for the next run, or else to change (slightly) one of the available parameters.
	
	\end{quote}

\smallskip
\noindent As a challenge for deep learning we suggest correspondingly:

\smallskip\begin{quote}
\noindent {\bf Challenge 1} \label{Chall1} Design a network which learns its sequential calibration for a strategy aiming at an average rank strictly below $U=2.33.$ \end{quote}

\smallskip 
Note that this is a suggested architecture \`a la  "play the winner rule." Calibration must thus probably be done for more than one parameter. Specialist in A.I. would possibly propose a more efficient approach. As we understand, deep learning does not really need such suggestions. (However, we do not see either why the assembling of black boxes, which create neural networks and for deep learning, could not profit from them.)

 An alternative challenge will be presented in Section \ref{Diffapp}. 
 
 \medskip
Allaart and Allen (2019), Assaf and Samuel-Cahn  (1996),  Christensen et al.(2013), Gnedin and Miretskiy (2007), Vanderbei (1980) and  many other people interested in Robbins' problem had intuitions which more or less into the right direction, and they had selected their objectives and adaptations accordingly.   

I have not heard or seen that notions equivalent to pre- or post-clouds had been advertized. However, Meier and S\"ogner (2017) were probably guided in their choice of a strategy by a similar idea. They obtained a simple non-trivial example combining rank dependent rules with threshold rules and attained an expected rank lower than the best upper bounds obtained in the literature so far, namely $2.32614....$

\subsection{Back to intuition}Applying specific strategies should not be able to compete with learning successive objectives and adaptations with the speed of a computer. Intuition tells us there should be more to it if we  argue as follows:

 Suppose we would like to find the minimum of a  complicated function $C: \R^m \to \R$ where each component of the solution vector in $\R^m$ is submitted to constraints resulting from its other components. We may have a "shot"
which hits the image not too far from the minimum. We succeeded in doing so
by beginning with an $ml$-strategy. Now look  at the parameters in Section 10.
 If deep learning
constantly updates parameters, it mimics what we could call "hyper-planes" of shots, producing a higher probability to cut the image of $C$  nearer to  the searched minimum.
\subsubsection {Why may deep learning show us patterns which our intuition is likely to overlook?}

Suppose $v<\tilde v.$ This means that the law of large numbers applied to the development of order statistics in increasing samples ${X_1, X_2, \cdots X_n, }$ (as $n$ grows), which led to the idea of $ml$-strategies, is too coarse to identify the limiting optimal strategy. This again implies that it either may be better to refuse an observation, although below the corresponding $ml$-threshold, or else be better to accept the same observation, although above the $ml$-threshold.
Our intuition grasps reasons why one or the other can be true. 
Again we should point out, that this can only 
work if $v < \tilde v$ B. and Ferguson had already in 1993 two (heuristic) arguments why they believe that $v<\tilde v$ and they never found good reason to change opinion.

Can deep learning give an answer? My feeling is that if deep-learning finds strategies for which for growing n the values $v_n$ stay clearly below the mentioned upper bound $U$ then a major step is done, because then $v <\tilde  v$ becomes increasingly probable. If not, even then this would be of interest, because this would tell us  that the law of large numbers takes over because the influence of the history fades away.

\subsection{Identifiability}

If deep-learning suggests an answer, then one must still confront an identifiability problem. Namely, one would like to decipher what the machine actually does to get the better values, i.e. one would wish to describe in words the essence of what modification of the $ml$-strategy has produced the improvements. 
This may still be difficult. We have only modest ideas how to tackle this problem, such as working again with the matrices of compatibility referred to in the Introduction (subsection \ref{compatibility}). A.I. specialists may then know how to apply a translater
into human speech. Note that our suggested approach \ref{Chall1} was guided by intuition. This is in contrast to the following alternative.

\section{An approach via a differential equation}
\label{Diffapp}
When Ben Green and Terence Tao received the Ostrowski Price at the Dutch Mathematical Congress (2007, U. Leiden), 
the author had the honour and pleasure to briefly discuss Robbins' problem with Professor Tao. 
Tao made two suggestions: 

\begin{quote}

{\it T1. Give the problem more structure.}
 
 \smallskip
 {\it T2. Try to find an equivalent problem.} 
 \end{quote}

I  probably failed to understand the essence of suggestion {\it T1},
and I still wonder how to give more structure to a problem with a clear definition without changing the problem as such. 
The second suggestion {\it T2} was  clear and encouraging.  Indeed, B. and Swan had, independently, already been trying to do this. We were looking at Robbins' problem in a continuous time embedding. This was more work than expected and required fine analytical steps, but  finally it led indeed to a fully equivalent problem (B. and Swan (2009)).
We present an outline of it since it yields an interesting alternative  to attack Robbins' problem. 
The setting is as follows:

\smallskip
{\narrower {\it Continuous time setting of Robbins' Problem:}  A decision maker observes  the $X_k\in[0,1]$ according to a planar Poisson process $(T_1, X_1),$ \,$ (T_2,X_2), \cdots$ with homogeneous rate $1$ on  ${\bf R}^+\times [0,1]$ and at most one of the $X_k$ can be accepted until time $t>0$. Acceptance is only possible at an arrival time $T_k$. The corresponding loss of an accepted $X_k$ is its absolute rank among those arriving up to time $t$. If no decision is taken up to time $t$, then the loss is given by some nonnegative  function $\Pi(t).$ What decision rule minimizes the expected loss? \par}

\smallskip
We first describe the set of admissible stopping rules. (See also Fathan  and Delage  (2021) for the notion of deep reinforcement learning for optimal stopping.)

Let $N(u)_{u\ge 0}$
denote the counting process of arrivals.  At each arrival time $T_{k}$ we can then use  all information up to time $T_k$. Formally this means, a stopping rule $\tau$ must satisfy $\{\tau \le s\} \in \F_s$ where $\s\{ ., ., ....,.\}$ denotes the sigma-field generated by the respective random variables listed in the argument. Moreover, the definition 
\begin{align}\F_s=\s\{(N_u)_{0\le u\le s}, (T_1,X_1), \cdots, (T_N(s), X_{N(s)})\} \end{align} is supposed to hold
with the understanding that $\F_s=\s\{(N_u)_{0\le u\le s}\}$ if there were no
arrivals up to time $s.$  (K\"uhne und R\"uschendorf (2016) and Gnedin (2007) called such stopping rules {\it canonical} stopping rules.) 

The value function  is thus accordingly
\begin{align} w(t)=\inf_\tau \E\left(R^{(t)}{\bf 1}_{\{T_\tau\le t\}}+\Pi(t){\bf 1}_{\{T_\tau>t\}}\right),
\end{align}
where $\Pi(t)$ is a penalty for not having stopped by time $t.$ (The total number of values arriving up to time $t$ is now a random variable denoted $N(t)$.) We do not specify this penalty function, but suppose that $\Pi(0)=0$ and that $\Pi(t)$ is Lipschitz continuous in $t$.

The continuous time counterpart of a memoryless threshold strategy as defined in B. and Ferguson (1993) is now the rule
\begin{align}\tau=\inf\{i \ge 1: X_i \le \varphi(T_i)\}
\end{align}
where $\varphi$ is  a real-valued positive function not depending on the history up to time $t.$ In analogy to B. and Ferguson (1993) we focus on a threshold function of the form  \begin{align}\label{varphi}\varphi(t) = 1_{\{0\le s\le t\}}\frac{c}{(1-s+c) }+1_{\{s>t\}}.\end{align}  Then from the planar Poisson process assumption we have $$P(T_\tau>s) =e^{-\mu(s)}$$
where $\mu(s)=\int_0^s\varphi(u) du$ for $s\le t$.
As shown in B. and Swan (2009, Section 2.1), we can then write down the value function $W_\tau(t)$
yielding
\begin{align} W_\tau(t)=1+(\Pi(t)-1)e^{-\mu(t)}+\frac{1}{2}\int_0^t\varphi(s)^2(t-s)e^{-\mu(s)}ds\nonumber\\+\frac{1}{2}\int_0^t\int_0^s\frac{(\varphi(s)-\varphi(u))^2}{1-\varphi(u)}du \,e^{-\mu(s)}ds~~~~~~~~\nonumber\end{align}
A finer analysis shows finally that the value function $w(t):=W_{\tau_t^*}$ under the optimal stopping rule $\tau_t^*$ with respect to the horizon $t$ is a differentiable function of $t$ and satisfies \begin{align}w'(t)+w(t)=\int_0^1\min\{1+xt,w(t|x)\} dx +\chi(t) \end{align}
where now $w(t|x)$ denotes the optimal value conditioned on a first value $x$ observed at time $0$ {\it which cannot be accepted}. A further analysis shows that $\chi(t)\to 0$ as $t\to\infty$ so that, by writing $h(t,x)=w(t|x)-w(t)$ with the constraint $h(t,x)\to 0$ as $t \to \infty$ we can limit our interest to
\begin{align}\label{Diff} w'(t)+w(t)=\int_0^1\min\{1+xt, w(t)+h(t,x)\} dx.\end{align}
This is a differential equation in the two unknown functions $w(t)$ and $h(t,x)$ and the challenge is to find a good estimate for $h(t,x).$  This brings us to the alternative challenge.
\medskip
\begin{quote}
{\bf Challenge 2}  \label{Chall2} Design an deep-learning approach  estimating $h(t,x)$ and then solve (\ref{Diff}) to obtain $\lim_{t\to\infty} w(t)).$
\end{quote}

\smallskip
The advantage of this closed form equation (\ref{Diff}) is that parametrization and calibration evoked in Challenge 1 is now implicitly resumed in one estimation problem. This is a bonus of Challenge 2. A malus may be the fact that our intuition seems lost in the differential equation (\ref{Diff}) with two unknown functions.

\section{General interest of full-history-dependent problems}Robbins' problem exemplifies a whole class of hard problems. The main part of its difficulty stems from the full-history dependence of the optimal strategy, and this can be a true problem in contexts other than that of selection strategies.

Learning is often based on huge data sets stemming from patients and certain control sets. Search can be accelerated if you know better at what patterns one should look first, and statistical inference becomes important. When continuous variables may trigger off discontinuous, things may become difficult.
In Robbins' problem, even the slightest change in the position of a point can change the loss by 1 (since ranks are integers.) The problem can be similar when recording counting processes. 

For instance, when searching for new drugs, in recording data (of a previously applied medication)  one may observe a dosage just slightly above the threshold necessary to kill all bacteria of a certain kind whereas a slightly insufficient dosage may give bacteria the chance to survive.  This may bias the success count.

\section {Why this article?}
 How come that an author who is a layman in deep learning, and who does not know too much about Robbins' problem either, writes an article on both subjects.

\smallskip
The answer is that the author sees a chance that a skilful deep learning input will solve the import question "$v= \tilde v?$",  since the main problem is  a limit-problem {\it without} a  $\forall$-statement, and since it is {\it easy to generate arbitrarily large training sets} of i.i.d. uniform random variables. 
As  shown in Section 7 (A) and (B), the question  is unlikely to be solved directly by other methods.
If deep learning, implemented  in Challenge 1 or Challenge 2, yields
a limiting average accepted rank smaller than $U,$ then this would mean real progress on Robbins' problem, and the value $v$ may be obtainable over short or long with improved deep learning. 

If not, then deep learning would also yield valuable information. Indeed, if the pre-cloud/post cloud arguments show no effect under better and better deep-learning, it becomes increasingly likely that $v=\tilde v.$ Keeping in mind the strong correlation between values and ranks shown in  (\ref{corr} ),
this would mean for probabilists a heart-warming conclusion : It is that with high probability, depending on the depth of deep learning, our intuition did not fail us. In other words, with high probability the strong law of large numbers combined with (\ref{corr}) has a  long arm, and we can confine our interest to $ml$-strategies.

\section{Conclusion}

Deep-learning is expected to be more than a {\it d\'ej\`a vu}. It would be nice to see this exemplified by progress on Robbins' problem which has caught the attention of authors in the domain of optimal stopping since 1990. 

 In our Introduction we have also discussed the question whether, true or false, computers are important, but not for Mathematics. (Of course, the importance for numerical computation has always been clear. )We have gone through several stages to try to answer but we never gave one. Could we do so?
 
If intuition allows us, with the help of computers, to push intuitions on 
 a new level in which we can check new hypotheses  (as e.g. our arguments based on pre-clouds and post-clouds in (\ref{prep1}) , (\ref{prep2}) and  (\ref{prep3}), and this with the parameters in  \ref{Chall1}  we want to choose, then I  would personally not question the importance of computers in Mathematics. 
Since we do not know before we look at a problem, whether a computer may have this potential,  we cannot say more. This is no cowardice. Different people, depending on the kind of problems they are studying, will have, rightly, different answers. 

There is also something joyful in our conclusion. We can afford to have different answers because we have different perceptions where real Mathematics begins, we have different intuitions, and, in particular,  so  much freedom in choosing problems we consider worth studying.

\section*{Acknowledgement} The author thanks G. Bontempi and J.-F. Raskin for discussions on deep learning, and  M. Duerinckx and P. Ernst for earlier discussions on Robbins' problem.

\section*{References}

~~~Acemoglu D. and Restrepo P (2019), {\it Artificial Intelligence, Automation and Work}, in
The Economics of Artificial Intelligence: An Agenda, edited by A. Agrawal, J. Gans and A. Goldfarb, University of Chicago Press, 197-236.


\smallskip 
Allaart P.C and Allen A. (2019)) {\it A random walk version of Robbins' problem: small horizon}, Math. Applicanda, Vol. 47, No. 2, 293-312 (2019)

\smallskip
Appel K. and Haken W. (1977), {\it Every Planar Map is Four Colorable. I. Discharging}, Illinois Journal of Mathematics, Vol. 21 (3): 429-490.

\smallskip
Assaf and Samuel-Cahn  (1996), {\it The secretary problem: minimizing the expected rank with I.I.D. random variables},  Adv. Appl. Prob.,Vol. 28, 828 - 852. 

\smallskip 
Bay\'on L., Fortuny Ayuso P.,  Grau Ribas J. M., Oller-Marc\'en A.M. and Ruiz M. M. (2018),
{\it The Best-or-Worst and the Postdoc problem},
J. Combinatorial Optimization 35, 703-723.

\smallskip Bruss F.T. (1982)  {\it A probabilistic approach to an approximation problem}, Ann. de la Soc. Scient. de Bruxelles, Vol. 96, Ser. I, 91-97.

\smallskip Bruss F.T. (1984)  {\it A unified approach to  a class of best choice problems with an unknown number of options}, Ann. of Probab.,  Vol. 12, 882-889.

\smallskip Bruss F.T. (1992) {\it Optimal Selection, Learning and Machine Implementation}), Contemporary Mathematics, Amer. Math. Society, Vol. 125, 3-35.

\smallskip Bruss F.T. (2005) {\it What is known about Robbins' problem?}, J. Appl. Prob., Vol. 42, 108-120.

\smallskip Bruss F.T. and Samuels S.M. (1987), {\it A unified approach to  a class of optimal selection problems with an unknown number of options}, Ann. of Probab.,  Vol. 15, 824-830.

\smallskip Bruss F.T. and Samuels S.M. (1990) {\it Conditions for quasi-stationarity of the Bayes rule in selection problems with an unknown number of rankable options}, Ann. of Probab.,  Vol. 18, 877-886.

\smallskip Bruss F.T. and Ferguson T.S. (1993), {\it Minimizing the expected rank with full information}, J. Appl. Prob.,Vol. 30, 616-626. 

\smallskip Bruss F.T. and Ferguson T.S. (1996), {\it Half-prophets and Robbins' problem of minimizing the expected rank}, (Athens Conference on Applied Probability and Time Series), Springer  Lecture Notes in Statistics, Vol.114 (I), 1-17.

\smallskip
Bruss F.T. and  Swan Y.C.  (2009) {\it The Poisson-embedded version of Robbins' problem of minimizing the expected rank.}, J. Appl. Prob., Vol. 46, 1-18.

\smallskip Chow Y.S, Moriguti S., Robbins H. and Samuels S.M. (1964), {\it Optimal selection based on relative rank (the secretary problem)}, Israel J. Math., Vol. 2, 81-90.

\smallskip Courant R. and Robbins H. (1941), {\it What is Mathematics?} Oxford University Press, (1st edition).

\smallskip Christensen S., Salminen P. and Ta B.Q. (2013) {\it Optimal stopping of strong Markov processes}, Stoch. Proc. Applic., Vol. 123 (3), 1138-1159.

\smallskip De la Vall\'ee Poussin A.  (1897), {\it Th\'eorie analytique des nombres premiers}, Ann. de la Soc. Scient. de Bruxelles, Vol. 20, Part II, 153- 250,  281-397.

\smallskip
Dendievel R. and Swan Y.C. (2016), {\it One step more in Robbins' Problem: Explicit solution for the case $n=4$}, Math. Applicanda,
Vol. 44  No 1, 135-148.

\smallskip
Fathan A. and Delage E. (2021), {\it Deep reinforcement learning for optimal stopping with application in financial engineering}, arXiv:2105.08877v1.


 
\smallskip
Gnedin A.  (2007) {\it Optimal Stopping with Rank-Dependent Loss}, J. Appl. Prob. Vol. 44, 996-1011.


\smallskip
Gnedin A. and Miretskiy D. (2007)   {\it Winning rate in the full-information best-choice problem},  J. Appl. Prob. Vol. 44,
560-565.

\smallskip 
Gnedin A. and Iksanov A. (2011) {\it Moments of Random Sums and Robbins' Problem of Optimal Stopping}, J. Appl. Prob., Vol. 48, 1197 - 1199. 

\smallskip Goldenshluger A. ,  Malinovsky Y. , and Zeevi A. (2020),
{\it A unified approach for solving sequential selection problems},
Probab. Surveys 17, 214-256.

\smallskip
Goodfellow I., Bengio Y. , Courville A.  (2016), {\it Deep Learning}, MIT Press.

\smallskip Grau Ribas J.M.  (2020)
An extension of the Last-Success-Problem,
Statistics \& Probability Letters,
Volume 156.

\smallskip K\"uhne R. and R\"uschendorf L. (2000), {\it Approximation of optimal stopping problems}, Stoch. Proc,. Applic.,
Vol. 90, 303-325.

\smallskip K\"uhne R. and R\"uschendorf L. (2003), {\it Approximate optimal stopping of dependent sequences}, Theory Probab. Appl., Vol. 48 (3), 465-480.

\smallskip
Li X. and Lee C.G.  (2023) {\it V-learning: An adaptive reinforcement learning algorithm for the optimal stopping problem},
Expert Systems with Applications,
Acc. 3 June 2023, Avail. online 8 June 2023.

\smallskip
Meier M. and S\"ogner L.  (2017),
{\it A new strategy for Robbins' problem of optimal stopping}, J. Appl. Prob., Vol. 54, 331 - 336.

\smallskip Nicholson, C. (2020) {\it A beginners guide to deep learning and neural networks},
https://wiki.pathmind.com/neural-network

\smallskip
Ramsey D.M. (2016) {\it A secretary problem with missing observations}, 
Math. Applicanda, Vol. 44 (1), 149-166.

\smallskip
Rutten X. (1987) {\it Etude et optimisation des fonctions  d'utilit\'e g\'en\'erales.}, M\'emoire de licence en sciences math\'ematiques, Facult\'e des sciences, Universit\'e de  Namur.

\smallskip 
R\"uschendorf L. (2016) {\it Approximative solutions of optimal stopping and selection problems}, Math. Applicanda, Vol. 44 (1), 17-44.

\smallskip
Simonis E. (1987) {\it Probl\`emes de s\'election optimale: Strat\'egies Autodidactes}, M\'emoire de licence en sciences math\'ematiques,  Facult\'e des sciences, Universit\'e de  Namur.

\smallskip
Szajowski K. (2007) {\it A game version of the Cowan-Zabczyk-Bruss problem}, Stat.\& Probab. Letters, Vol 77, 1683-1689.

\smallskip 
Szajowski K. (2009) {A rank-based selection with cardinal payoffs and a cost of choice}, Sc. Mat. Jap. Vol 69 (2), 69-77. 

\smallskip
Vanderbei R.J. (1980) {\it Optimal choice of a subset of a population.}, Math. Oper. Res., Vol 5, 481-486.

\smallskip
Villeneuve S. (2007) {\it On threshold strategies and the smooth-fit principle for optimal stopping problems.} J. Appl. Prob., Vol. 44, 181-198. 

\smallskip
Wang X., Han Y., Leung VCM, Niato D., Yan X. and Chen X. (2020), {\it Convergence of edge computing and deep learning: A comprehensive survey}, IEEE Communications Surveys \& Tutorials, Vol. 22, 869-904.

\bigskip

\noindent F. Thomas Bruss\\
Universit\'e Libre de Bruxelles\\
Department of Mathematics\\
Campus Plaine, CP 210\\
B-1050 Brussels, Belgium

\end{document}